\documentstyle[12pt, amsfonts]{amsart}
\begin{document}
\def\R{{\mathbb R}}
\def\Z{{\mathbb Z}}
\def\C{{\mathbb C}}
\newcommand{\trace}{\rm trace}
\newcommand{\Ex}{{\mathbb{E}}}
\newcommand{\Prob}{{\mathbb{P}}}
\newcommand{\E}{{\cal E}}
\newcommand{\F}{{\cal F}}
\newtheorem{df}{Definition}
\newtheorem{theorem}{Theorem}
\newtheorem{lemma}{Lemma}
\newtheorem{pr}{Proposition}
\newtheorem{co}{Corollary}
\def\n{\nu}
\def\sign{\mbox{ sign }}
\def\a{\alpha}
\def\N{{\mathbb N}}
\def\A{{\cal A}}
\def\L{{\cal L}}
\def\X{{\cal X}}
\def\F{{\cal F}}
\def\c{\bar{c}}
\def\v{\nu}
\def\d{\delta}
\def\diam{\mbox{\rm dim}}
\def\vol{\mbox{\rm Vol}}
\def\b{\beta}
\def\t{\theta}
\def\l{\lambda}
\def\e{\varepsilon}
\def\colon{{:}\;}
\def\pf{\noindent {\bf Proof :  \  }}
\def\endpf{ \begin{flushright}
$ \Box $ \\
\end{flushright}}
%%%%%%%%%%%%%%%%%%%%%%%%%%%%%%%%%%%%%%%%%%%%%%%%%%%%%%%%%%%%%%%%%%%%

\title[Stability and separation]{Stability and separation in volume comparison problems}

\author{Alexander Koldobsky}

\address{Department of Mathematics\\ 
University of Missouri\\
Columbia, MO 65211}

\email{koldobskiya@@missouri.edu}

%%%%%%%%%%%%%%%%%%%%%%%%%%%%%%%%%%%%%%%%%%%%%%%%%%%%%%%%%%%%%%%%%%%%
\begin{abstract}   We review recent stability and separation results in volume comparison problems
and use them to prove several hyperplane inequalities for intersection and projection bodies.
\end{abstract}  
\maketitle
%%%%%%%%%%%%%%%%%%%%%%%%%%%%%%%%%%%%%%%%%%%%%%%%%%%%%%%%%%%%%%%%%%%%

\section{Introduction}

A typical {\it volume comparison problem} 
asks whether inequalities $$f_K(\xi)\le f_L(\xi), \qquad  \forall \xi\in S^{n-1}$$ imply $|K|\le |L|$ 
for any $K,L$ from a certain class of origin-symmetric convex bodies in $\R^n,$ where $f_K$ 
is a geometric characteristic of $K.$ One can have in mind the hyperplane section function
$f_K(\xi)= |K\cap \xi^\bot|,$ where $|K|$ stands for volume of proper dimension and 
$\xi^\bot$ is the central hyperplane perpendicular to $\xi\in S^{n-1}.$

In the case where the answer to a volume comparison problem is affirmative, one can ask
a stronger {\it stability} question. Suppose that $\e>0$ and 
\begin{equation}\label{stab1}
f_K(\xi)\le f_L(\xi)+\e, \qquad  \forall \xi\in S^{n-1}.
\end{equation}
Does there exist a constant $c$ not dependent on $\e$ and such that for every $\e>0$
\begin{equation}\label{volcomp1}
|K|^{\frac{n-1}n} \le |L|^{\frac{n-1}n} + c\e ?
\end{equation}

Stability results are related to hyperplane inequalities as follows. Suppose stability holds 
for both pairs $K,L$ and $L,K$ with the same constant $c.$ Put 
$$\e= \max_{\xi\in S^{n-1}} \left|f_K(\xi)-f_L(\xi) \right|,$$
then one can switch $K$ and $L$ in (\ref{stab1}) and, correspondingly, in (\ref{volcomp1}).
The resulting inequality for volumes will be called a {\it volume difference inequality}:
\begin{equation}\label{diff1}
\left| |K|^{\frac{n-1}n} - |L|^{\frac{n-1}n}\right| \le c\e = c \max_{\xi\in S^{n-1}} \left|f_K(\xi)-f_L(\xi) \right|.
\end{equation}
Suppose now that the function $f_L$ converges to zero uniformly with respect to $\xi$ when $L=\beta B_2^n$ is
a multiple of the unit Euclidean ball and $\beta\to 0.$  Then, when $L=\beta B_2^n$ and $\beta \to 0,$ the inequality (\ref{diff1}) turns into 
what we call a {\it hyperplane inequality}:
\begin{equation}\label{hyper1}
|K|^{\frac{n-1}n} \le c\max_{\xi\in S^{n-1}} f_K(\xi).
\end{equation}
\medbreak
One can also consider a {\it separation} problem. Suppose that $\e>0$ and 
\begin{equation}\label{stab2}
f_K(\xi)\le f_L(\xi)-\e, \qquad  \forall \xi\in S^{n-1}.
\end{equation}
Does there exist a constant $c$ not dependent on $\e$ and such that for every $\e>0$
\begin{equation}\label{volcomp2}
|K|^{\frac{n-1}n} \le |L|^{\frac{n-1}n} - c\e ?
\end{equation}
In the case where the answer is affirmative, assuming that
$$\e= \min_{\xi\in S^{n-1}} \left(f_L(\xi) - f_K(\xi)\right)>0,$$
we get another kind of a volume difference inequality:
\begin{equation}\label{diff2}
 |L|^{\frac{n-1}n} - |K|^{\frac{n-1}n} \ge c\e = c \min_{\xi\in S^{n-1}} \left(f_L(\xi)-f_K(\xi)\right).
\end{equation}
Again, if $f_{\beta B_2^n}$ converges to zero uniformly in $\xi$ when $\beta \to 0,$ 
we get the following version of  a hyperplane inequality:
\begin{equation}\label{hyper2}
|L|^{\frac{n-1}n} \ge c\min_{\xi\in S^{n-1}} f_L(\xi).
\end{equation}

This strategy was first applied in \cite{K6} to several functions $f_K$ including the hyperplane section function 
and the hyperplane projection function.
In \cite{K8} similar inequalities were proved for arbitrary measure with continuous density in place of volume.
Sections of lower dimensions were considered in \cite{KM}, and stability and hyperplane
inequalities for complex convex bodies were proved in \cite{K7, KPZ}. 

In this article we review stability and separation results and prove some of them with the best possible
constants,  while in the original papers the constants were sometimes estimated. The proofs are based on recently
developed Fourier analytic approach to sections and projections of convex bodies; see \cite{K4, KRZ, KY}.
We also prove several hyperplane inequalities for intersection and projection bodies.

\section{Hyperplane sections}

Suppose that 
$$f_K(\xi)=S_K(\xi)= \left|K\cap\xi^\perp\right|, \qquad \xi\in S^{n-1},$$
is the hyperplane section function, then the volume comparison question 
is the matter of the Busemann-Petty problem, raised in 1956 in \cite{BP}.
Let $K,L$ be origin-symmetric convex bodies in $\R^n$ such that 
$\left|K\cap\xi^\perp\right| \le \left|L\cap\xi^\perp\right|$ for 
every $\xi\in S^{n-1}.$ 
Does it necessarily follow that $\left|K\right| \le \left|L\right| ?$
The problem was solved at the end of the 1990's as the result 
of a sequence of papers \cite{LR}, \cite{Ba1}, \cite{Gi}, \cite{Bo4}, 
\cite{L}, \cite{Pa}, \cite{G1}, \cite{G2}, \cite{Z1}, \cite{Z2}, \cite{K2}, \cite{K3}, \cite{Z3},
\cite{GKS} ; see \cite[p. 3]{K4} or \cite[p. 343]{G3} for the history of the solution.
The answer is affirmative if $n\le 4$, and it is negative if $n\ge 5.$ Moreover,
Lutwak \cite{L} proved that if $K$ is an {\it intersection body} (see definition below) and $L$ is any origin-symmetric
star body, then the answer to the Busemann-Petty problem is affirmative in every dimension.

The corresponding stability result was proved in \cite[Theorem 1]{K6}. The theorem is stated in \cite{K6}
with $c_n$ replaced by 1, though the proof there actually establishes the result with the constant $c_n,$ which is the 
best possible.  Also, the proof in \cite{K6} is geometric, while here we use methods of Fourier analysis.

Throughout the paper 
$$c_n:= \frac{|B_2^n|^{\frac {n-1}n}}{|B_2^{n-1}|},$$
where $B_2^n$ is the unit Euclidean ball. Note that  $c_n\in (\frac 1{\sqrt{e}},1);$ see for example \cite[Lemma 2.1]{KL}.

\begin{theorem} {\bf (\cite{K6})} \label{main-int} Suppose that $\e>0$,  $K$ and $L$ are origin-symmetric
star bodies in $\R^n,$ and $K$ is an intersection body.  If for every $\xi\in S^{n-1}$
\begin{eqnarray}\label{sect1}
|K\cap \xi^\bot| \le  |L\cap \xi^\bot| + \e,
\end{eqnarray}
then
$$|K|^{\frac{n-1}n}  \le |L|^{\frac{n-1}n} + c_n \e.$$
\end{theorem}
Recall that the constant $c_n<1.$
To prove Theorem \ref{main-int} we need several definitions and known facts.
We say that a closed bounded set $K$ in $\R^n$ is a {\it star body}  if 
every straight line passing through the origin crosses the boundary of $K$ 
at exactly two points different from the origin, the origin is an interior point of $K,$
and the {\it Minkowski functional} 
of $K$ defined by 
$$\|x\|_K = \min\{a\ge 0:\ x\in aK\}$$
is a continuous function on $\R^n.$ 

The {\it radial function} of a star body $K$ is defined by
$$\rho_K(x) = \|x\|_K^{-1}, \qquad x\in \R^n.$$
If $x\in S^{n-1}$ then $\rho_K(x)$ is the radius of $K$ in the
direction of $x.$

Writing volume in polar coordinates we get the polar formula for volume
\begin{equation} \label{polar-volume}
|K|
=\frac{1}{n} \int_{S^{n-1}} \rho_K^n(\theta) d\theta=
\frac{1}{n} \int_{S^{n-1}} \|\theta\|_K^{-n} d\theta,
\end{equation}
and the polar formula for the volume of a section
\begin{equation}\label{polar-section}
|K\cap \xi^\bot| = \frac 1{n-1} \int_{S^{n-1}\cap \xi^\bot} \rho_K^{n-1}(\theta) d\theta.
\end{equation}

The class of {\it intersection bodies} was introduced by Lutwak \cite{L}.
Let $K, L$ be origin-symmetric star bodies in $\R^n.$ We say that $K$ is the 
intersection body of $L$ if the radius of $K$ in every direction is 
equal to the $(n-1)$-dimensional volume of the section of $L$ by the central
hyperplane orthogonal to this direction, i.e. for every $\xi\in S^{n-1},$
\begin{equation} \label{intbodyofstar}
\rho_K(\xi)= \|\xi\|_K^{-1} = \left|L\cap \xi^\bot\right|.
\end{equation} \index{intersection body of a star body}
A more general class of {\it intersection bodies} can be defined (see \cite{GLW}) as the closure of the
class of intersection bodes of star bodies in the radial metric 
$$\rho(K,L)= \max_{\xi\in S^{n-1}} \left| \rho_K(\xi)-\rho_L(\xi)\right|.$$

We consider Schwartz distributions, i.e. continuous functionals on the space ${\cal{S}}(\R^n)$
of rapidly decreasing infinitely differentiable functions on $\R^n$. 
The Fourier transform of a distribution $f$ is defined by $\langle\hat{f}, \phi\rangle= \langle f, \hat{\phi} \rangle$ for
every test function $\phi \in {\cal{S}}(\R^n).$ For any even distribution $f$, we have $(\hat{f})^\wedge
= (2\pi)^n f$.

If $K$ is a star body  and $0<p<n,$
then $\|\cdot\|_K^{-p}$  is a locally integrable function on $\R^n$ and represents a distribution acting by integration. 
Suppose that $K$ is infinitely smooth, i.e. $\|\cdot\|_K\in C^\infty(S^{n-1})$ is an infinitely differentiable 
function on the sphere. Then by \cite[Lemma 3.16]{K4}, the Fourier transform of $\|\cdot\|_K^{-p}$  
is an extension of some function $g\in C^\infty(S^{n-1})$ to a homogeneous function of degree
$-n+p$ on $\R^n.$ When we write $\left(\|\cdot\|_K^{-p}\right)^\wedge(\xi),$ we mean $g(\xi),\ \xi \in S^{n-1}.$
If $K,L$ are infinitely smooth star bodies, the following spherical version of Parseval's
formula was proved in \cite{K5} (see \cite[Lemma 3.22]{K4}):  for any $p\in (-n,0)$
\begin{equation}\label{parseval}
\int_{S^{n-1}} \left(\|\cdot\|_K^{-p}\right)^\wedge(\xi) \left(\|\cdot\|_L^{-n+p}\right)^\wedge(\xi) =
(2\pi)^n \int_{S^{n-1}} \|x\|_K^{-p} \|x\|_L^{-n+p}\ dx.
\end{equation}

A distribution is called {\it positive definite} if its Fourier transform is a positive distribution in
the sense that $\langle \hat{f},\phi \rangle \ge 0$ for every non-negative test function $\phi.$
It was proved in \cite[Theorem 1]{K2} that an origin-symmetric star body in $\R^n$ is an
intersection body if and only if the function $\|\cdot\|_K^{-1}$ represents a positive definite 
distribution. As proved in \cite{G2,Z3} (see also \cite{GKS} or \cite[p. 73]{K4}), every origin-symmetric convex body in 
$\R^n, \ n\le 4$ is an intersection body. It was shown in \cite[Theorem 3]{K2} that the unit ball of 
any finite dimensional subspace of $L_p,\ 0<p\le 2$ is an intersection body. For other results
on intersection bodies, see \cite[Chapter 8]{G3} and \cite[Chapter 4]{K4}.

For origin-symmetric star bodies $K,L$ in $\R^n,$  the radial sum  $K+_r L$  of $K$ and $L$ 
is a star body defined by
$$\rho_{K+_rL}(\xi) = \rho_K(\xi) + \rho_L(\xi), \qquad \forall \xi\in S^{n-1}.$$
If $K$ and $L$ are both intersection bodies, then their radial sum is also an intersection body,
which follows, for example, from the Fourier characterization of intersection bodies formulated above.

\medbreak
\noindent {\bf Proof of Theorem \ref{main-int}.}  By approximation (see, for example \cite[Theorem 3.3.1]{S2}),
we can assume that the bodies $K$ and $L$ are infinitely smooth. It was proved in \cite{K1} that
\begin{eqnarray}\label{eqn:defS}
\left|K\cap \xi^\bot\right|=\frac{1}{\pi(n-1)}(\|x\|_K^{-n+1})^\wedge(\xi),\qquad \forall \xi\in S^{n-1},
\end{eqnarray}
so (\ref{sect1}) can be written as
\begin{equation} \label{form1}
(\|x\|_K^{-n+1})^\wedge(\xi) \le (\|x\|_L^{-n+1})^\wedge(\xi) + \pi(n-1)\e,\qquad  \forall \xi\in S^{n-1}.
\end{equation}
Also, by the remark before the proof and the Fourier characterization of intersection bodies,  
$(\|x\|_K^{-1})^\wedge$ is an infinitely smooth 
non-negative function on the sphere. By (\ref{form1}), the polar formula for volume (\ref{polar-volume}) 
and Parseval's formula on the sphere (\ref{parseval}),

$$(2\pi)^n n |K| = (2\pi)^n \int_{S^{n-1}}\|x\|_K^{-n+1}\|x\|_K^{-1}dx$$
$$=\int_{S^{n-1}}(\|x\|_K^{-1})^\wedge(\theta)(\|x\|_K^{-n+1})^\wedge(\theta)d\theta$$
$$\le \int_{S^{n-1}}(\|x\|_K^{-1})^\wedge(\theta)(\|x\|_L^{-n+1})^\wedge(\theta)d\theta$$
\begin{equation} \label{eq21}
+ {\pi(n-1)}\e \int_{S^{n-1}}(\|x\|_K^{-1})^\wedge(\theta)d\theta.
\end{equation}
By Parseval's formula and H\"older's inequality,
$$\int_{S^{n-1}}(\|x\|_K^{-1})^\wedge(\theta)(\|x\|_L^{-n+1})^\wedge(\theta)d\theta$$
\begin{equation}\label{eq22}
= (2\pi)^n \int_{S^{n-1}}\|x\|_L^{-n+1}\|x\|_K^{-1}dx \le (2\pi)^n n |K|^{\frac1n} |L|^{\frac{n-1}n}.
\end{equation}

To estimate the second summand in (\ref{eq21}), we use the formula for the Fourier transform 
(in the sense of distributions; see \cite[p.194]{GS})
$$\left(|x|_2^{-n+1}\right)^\wedge(\theta) = \frac{2\pi^{\frac {n+1}2}}
{\Gamma(\frac{n-1}2)} |\theta|_2^{-1}.$$
Again using Parseval's formula and then H\"older's inequality,
$$\int_{S^{n-1}}(\|x\|_K^{-1})^\wedge(\theta)d\theta$$
$$ = 
\frac{\Gamma(\frac {n-1}2)}{2\pi^{\frac {n+1}2}}
\int_{S^{n-1}}(\|x\|_K^{-1})^\wedge(\theta)\left(|x|_2^{-n+1}\right)^\wedge(\theta)d\theta$$
$$= \frac{(2\pi)^n \Gamma(\frac {n-1}2)}{2\pi^{\frac {n+1}2}}
\int_{S^{n-1}} \|x\|_K^{-1}\ dx$$
$$ \le \frac{(2\pi)^n \Gamma(\frac {n-1}2)|S^{n-1}|^{\frac{n-1}n}}{2\pi^{\frac {n+1}2}}
\left(\int_{S^{n-1}} \|x\|_K^{-n}\ dx\right)^{\frac 1n}$$
$$= \frac{(2\pi)^n \Gamma(\frac {n-1}2)|S^{n-1}|^{\frac{n-1}n}}{2\pi^{\frac {n+1}2}} \left(n|K|\right)^{\frac 1n}$$
Combining this with (\ref{eq21}) and (\ref{eq22}), we get
$$(2\pi)^n n |K| \le (2\pi)^n n|K|^{\frac1n} |L|^{\frac{n-1}n} $$
$$+ \frac{(2\pi)^n\e \pi(n-1)n^{\frac 1n}  \Gamma(\frac{n-1}2)
\left|S^{n-1}\right|^{\frac {n-1}n}}{2\pi^{\frac {n+1}2}} |K|^{1/n}.$$
Now to represent the coefficient in the required form use
$$|S^{n-1}|=n|B_2^n| = \frac {2\pi^{\frac n2}}{\Gamma(\frac n2)}.$$
\endpf

Interchanging $K$ and $L$ in Theorem \ref{main-int}, we get the corresponding volume difference inequality.

\begin{co}  \label{n4} If $K$ and $L$ are intersection bodies in $\R^n$ (in particular, any origin-symmetric 
convex bodies in $\R^3$ or $\R^4$), then 
$$\left| |K|^{\frac{n-1}n} - |L|^{\frac{n-1}n}\right| \le c_n \max_{\xi\in S^{n-1}} \left| |K\cap \xi^\bot| - |L\cap \xi^\bot|\right|.$$
\end{co}
\bigbreak
If $L=\delta B_2^n$ in the latter inequality, then sending $\delta$ to zero
we get that for any intersection body $K$ in $\R^n$
and, in particular, any origin-symmetric convex body in $\R^3$ or $\R^4,$
\begin{equation}\label{hyper-inter}
|K|^{\frac {n-1}n} \le c_n \max_{\xi \in S^{n-1}} |K\cap \xi^\bot|.
\end{equation}
Inequality (\ref{hyper-inter}) also immediately follows from the affirmative 
part of the Busemann-Petty problem; see \cite[Theorem 9.4.11]{G3}. 
Note that (\ref{hyper-inter}) is a particular case of the well-known and still open 
Hyperplane Problem (see \cite{Bo1, Bo2, Ba3, MP}) 
which can be formulated as follows.
Does there exist  an absolute constant $C$ so that for any origin-symmetric convex body $K$ in $\R^n$
$$
|K|^{\frac {n-1}n} \le C \max_{\xi \in S^{n-1}} |K\cap \xi^\bot|.
$$
The best-to-date estimate $C\sim n^{1/4}$ belongs
to Klartag \cite{Kl}, who slightly improved the previous estimate of Bourgain \cite{Bo3}.

The volume difference inequality of Corollary \ref{n4} can also be used to prove  a hyperplane inequality
for the average volume of central hyperplane sections, which we denote by
$${\rm{as}}(K) = \frac 1{|S^{n-1}|} \int_{S^{n-1}} |K\cap \xi^\bot| d\xi.$$
 For any continuous function $h$ on $S^{n-1},$
 $$|S^{n-2}| \int_{S^{n-1}} h(x) dx = \int_{S^{n-1}}\left(\int_{S^{n-1}\cap \xi^\bot} h(x)dx \right)d\xi.$$
Using this and (\ref{polar-section}) ,
$$ {\rm{as}}(K) = \frac 1{(n-1)|S^{n-1}|} \int_{S^{n-1}} \left(\int_{S^{n-1}\cap \xi^\bot} \rho_K^{n-1}(\theta) d\theta\right) d\xi$$
\begin{equation} \label{as-polar}
= \frac {|S^{n-2}|}{(n-1)|S^{n-1}|} \int_{S^{n-1}}  \rho_K^{n-1}(\theta) d\theta.
\end{equation}

\begin{co} \label{hyper-aversect} If $K$ is an intersection body in $\R^n, \ n\ge 3,$ then
$${\rm{as}}(K) \le \frac {|B_2^{n-1}|}{ |B_2^{n-2}|  |B_2^n|^{\frac 1n}}    \max_{\xi\in S^{n-1}} {\rm{as}}(K\cap \xi^\bot)
\  |K|^{\frac 1n},$$
 with equality when $K=B_2^n.$
\end{co}

\pf  Since $K$ and $B_2^n$ are intersection bodies, for every $\e>0$ the radial sum $K+_r \e B_2^n$ is also an intersection 
body. By Corollary \ref{n4} applied to the bodies $K+_r \e B_2^n$ and $K,$ we get that for every $\e>0$
\begin{equation} \label{radialsum}
\frac{|K+_r \e B_2^n|^{\frac{n-1}n} - |K|^{\frac{n-1}n}}{\e} \le 
c_n \max_{\xi\in S^{n-1}}\frac {|(K\cap \xi^\bot)+_r \e B_2^{n-1}| - |K\cap \xi^\bot|}{\e}.
\end{equation}
By the polar formula for the volume (\ref{polar-volume}),
$$|K+_r \e B_2^n| = \frac 1n \int_{S^{n-1}} (\rho_K(\theta)+\e)^n d\theta,$$
so 
$$\lim_{\e\to 0} \frac{|K+_r \e B_2^n|^{\frac{n-1}n} - |K|^{\frac{n-1}n}}{\e} = \frac {n-1}n |K|^{-\frac 1n} \int_{S^{n-1}} \rho_K^{n-1}(\theta)$$
$$=\frac {n-1}n |K|^{-\frac 1n} \frac {(n-1)|S^{n-1}|}{|S^{n-2}|}\  {\rm{as}}(K).$$
Similarly, the limit of the right-hand side of (\ref{radialsum}), as $\e\to 0,$ is equal to
$$c_n \frac {(n-2) |S^{n-2}|}{|S^{n-3}|} \max_{\xi\in S^{n-1}} {\rm{as}}(K\cap \xi^\bot).$$
It is easily seen that the convergence of the quotient in the right-hand side of (\ref{radialsum}) is uniform with respect
to $\xi,$ as $\e\to 0,$ so one can switch the limit and maximum.
Sending $\e$ to 0 in (\ref{radialsum}) and using
$|S^{n-1}|= n |B_2^n|$ we get the result.
\endpf
\bigbreak
A separation result for hyperplane sections was proved in \cite[Theorem 2]{K6}. The constant $c$ in this
result does not depend on $\e,$ but depends on the dimension and on the normalized inradius of
$K:$
$$r(K) = \frac{\min_{\xi\in S^{n-1}} \rho_K(\xi)}{|K|^{1/n}}.$$

\begin{theorem} \label{main-int1} (\cite{K6}) Let $K$ and $L$ be origin-symmetric
star bodies in $\R^n$ and $\e>0.$ Assume that $K$ is an intersection body.  
If for every $\xi\in S^{n-1}$
\begin{eqnarray}\label{sect2}
|K\cap \xi^\bot| \le |L\cap \xi^\bot| - \e,
\end{eqnarray}
then
$$|K|^{\frac{n-1}n}  \le |L|^{\frac{n-1}n} -  \sqrt{\frac{2\pi}{n+1}}\ r(K) \e.$$
\end{theorem}
\bigbreak
Since the answer to the Busemann-Petty problem is negative in most dimensions,  
one may ask what information about the hyperplane section function $S_K$ does allow to compare the volumes 
in all dimensions. An answer to this question was given in \cite{KYY}:  
for two origin-symmetric infinitely smooth bodies $K,L$ in $\R^n$ and $\alpha\in [n-4,n-1)$ 
the inequalities
\begin{eqnarray}\label{eqn:condition1}
(-\Delta)^{\alpha/2} S_K(\xi)\le (-\Delta)^{\alpha/2} S_L(\xi), \qquad \forall \xi\in S^{n-1}
\end{eqnarray}
imply that $|K|\le |L|,$ while for $\alpha<n-4$ this is not necessarily
true.
Here $\Delta$ is the Laplace operator on $\R^n$, and the fractional powers of the
Laplacian are defined by
$$
(-\Delta)^{\alpha/2}f = \frac{1}{(2\pi)^n}( |x|_2^\alpha  \hat{f}(x))^\wedge,
$$
where  the Fourier transform is considered in the sense of distributions, $|x|_2$ stands for
the Euclidean norm in $\R^n$, and the function $S_K$ is extended in (\ref{eqn:condition1}) 
to a homogeneous function of degree -1 on the whole $\R^n.$ The corresponding stability 
result was proved in \cite[Theorem 3]{K6}.

\begin{theorem} \label{main-kyy} (\cite{K6}) Let $\e>0,\ \alpha \in [n-4,n-1)$,  and let $K$ and $L$ be origin-symmetric
infinitely smooth convex bodies in $\R^n$, $n\ge 4$, so that for every $\xi\in S^{n-1}$
\begin{eqnarray}\label{eqn:condition11}
(-\Delta)^{\alpha/2} S_K(\xi)\le (-\Delta)^{\alpha/2} S_L(\xi) +\e.
\end{eqnarray}
Then
$$|K|^{\frac{n-1}n}  \le |L|^{\frac{n-1}n} + c \e,$$
where 
$$c=c(\alpha,n) = \frac{\sqrt{\pi}(n-1)\Gamma(\frac{n-\alpha-1}2)}{2^{\alpha+\frac1n} n^{\frac{n-1}n}
\Gamma(\frac{\alpha+1}2) \left(\Gamma(\frac n2)\right)^{\frac{n-1}n}}.$$
\end{theorem}

A separation result was proved in \cite[Theorem 4]{K6}.

\begin{theorem} \label{main2-kyy} (\cite{K6}) Let $\e>0,\ \alpha \in [n-4,n-1)$,  $K$ and $L$ be origin-symmetric
infinitely smooth convex bodies in $\R^n$, $n\ge 4$, so that for every $\xi\in S^{n-1}$
\begin{eqnarray}\label{eqn:condition}
(-\Delta)^{\alpha/2} S_K(\xi)\le (-\Delta)^{\alpha/2} S_L(\xi) -\e.
\end{eqnarray}
Then
$$|K|^{\frac{n-1}n}  \le |L|^{\frac{n-1}n} - c \e,$$
where 
$$c= r(K)\frac{\pi(n-1)\Gamma(\frac{n-\alpha-1}2)}{n2^\alpha \Gamma(\frac{\alpha+1}2)\Gamma(\frac n2)}.$$
\end{theorem}

\section{Hyperplane projections}

Now we pass to the hyperplane projection function
$$f_K(\xi)= P_K(\xi)= |K\vert\xi^\perp|, $$
where $K\vert\xi^\perp$ is the orthogonal projection of $K$ to the hyperplane $\xi^\perp.$
The corresponding volume comparison result is known as Shephard's problem,
which was posed in 1964 in \cite{Sh} and solved soon after that by Petty \cite{Pe} and 
Schneider \cite{S1}, independently. Suppose that $K$ and $L$ are origin-symmetric convex bodies in 
$\R^n$ so that  $|K\vert \xi^\bot|\le |L\vert \xi^\bot|$ for every $\xi\in S^{n-1}.$ Does it follow that
$|K|\le |L|?$ The answer if affirmative only in dimension 2. Both solutions use the fact that
the answer to Shephard's problem is affirmative in every dimension 
under the additional assumption that $L$ is a projection body; see definition below.

In the case of projections the constant in the stability result depends on the body and dimension, while the 
constant in the separation result does not. One can say that separation 
is a more natural property for projections than stability, while for sections it is the other way around.
Therefore, we start with a separation result
which was proved in \cite[Theorem 6]{K6}.  The constant $c_n$ in \cite{K6}
was at the last moment estimated from below by $1/\sqrt{e}$, so we now formulate and prove the
result with the best possible constant.
\begin{theorem} \label{main-proj1} (\cite{K6}) Suppose that $\e>0$,  $K$ and $L$ are origin-symmetric
convex bodies in $\R^n,$ and $L$ is a projection body.  If for every $\xi\in S^{n-1}$
\begin{eqnarray}\label{proj2}
|K\vert \xi^\bot|\le |L\vert \xi^\bot| - \e,
\end{eqnarray}
then
$$|K|^{\frac{n-1}n}  \le |L|^{\frac{n-1}n} - c_n \e.$$
\end{theorem}

To prove Theorem \ref{main-proj1} we need several more definitions and results from convex geometry. 
We refer the reader to \cite{S2} for details.

The {\it support function} of a convex body $K$ in $\R^n$ is defined by
$$h_K(x) = \max_{\{\xi\in \R^n:\|\xi\|_K=1\}} (x,\xi),\quad x\in \R^n.$$ 
If $K$ is origin-symmetric, then $h_K$ is a norm on $\R^n.$

The {\it surface area measure} $S(K, \cdot)$ of a convex body $K$ in 
$\R^n$ is defined as follows. For every Borel set $E \subset S^{n-1},$ 
$S(K,E)$ is equal to Lebesgue measure of the part of the boundary of $K$
where normal vectors belong to $E.$ 
We usually consider bodies with absolutely continuous surface area measures.
A convex body $K$ is said to have the {\it curvature function} 
$$ f_K: S^{n-1} \to \R,$$
if its surface area measure $S(K, \cdot)$ is absolutely 
continuous with respect to Lebesgue measure $\sigma_{n-1}$ on 
$S^{n-1}$, and
$$
\frac{d S(K, \cdot)}{d \sigma_{n-1}}=f_K \in L_1(S^{n-1}),
$$
so $f_K$ is the density of $S(K,\cdot).$

By the approximation argument of \cite[Th. 3.3.1]{S2},
we may assume in the formulation of Shephard's problem that the bodies 
$K$ and $L$ are such that  their support functions $h_K,\ h_L$ are 
infinitely smooth functions on $\R^n\setminus \{0\}$.
Using \cite[Lemma 3.16]{K4}
we get in this case that
the Fourier transforms $\widehat{h_K},\ \widehat{h_L}$ are the
extensions of infinitely differentiable functions on the sphere
to homogeneous distributions on $\R^n$ of degree $-n-1.$
Moreover, by a similar approximation argument (see also \cite[Section 5]{GZ}),
we may assume that  our bodies have absolutely continuous surface area 
measures. Therefore, in the rest of this section, $K$ and $L$ are 
convex symmetric bodies with infinitely smooth support functions and absolutely 
continuous surface area measures.

The following version of Parseval's formula was proved in \cite{KRZ} (see also \cite[Lemma 8.8]{K4}):
\begin{equation} \label{pars-proj}
\int_{S^{n-1}} \widehat{h_K} (\xi) \widehat{f_L}(\xi)\ d\xi =
(2\pi)^n \int_{S^{n-1}} h_K(x) f_L(x)\ dx.
\end{equation}

The volume of a body can be expressed in terms of its support function and 
curvature function:
\begin{equation}\label{vol-proj}
|K| = \frac 1n \int_{S^{n-1}}h_K(x) f_K(x)\ dx.
\end{equation}

If $K$ and $L$ are two convex bodies in $\R^n$ the {\it mixed volume} $V_1(K,L)$
is equal to  
$$V_1(K,L)= \frac{1}{n} \lim_{\e\to +0}
\frac{|K+\epsilon L|- |K|}{\e}.$$
We use the following
first Minkowski inequality (see \cite[p.23]{K4}):  
for any convex bodies $K,L$ in $\R^n,$ 
\begin{equation} \label{firstmink}
V_1(K,L) \ge |K|^{\frac {n-1}n} |L|^{\frac 1n}.
\end{equation}
The mixed volume can also be expressed in terms of the support and
curvature functions:

\begin{equation}\label{mixvol-proj}
V_1(K,L) = \frac 1n \int_{S^{n-1}}h_L(x) f_K(x)\ dx.
\end{equation}

Let $K$ be an origin-symmetric convex body in $\R^n.$ The {\it
projection body} $\Pi K$ of $K$ is defined as an origin-symmetric convex 
body in $\R^n$ whose support function in every direction is equal to
the volume of the hyperplane projection of $K$ to this direction: 
for every $\theta\in S^{n-1},$ 
\begin{equation} \label{def:proj}
h_{\Pi K}(\theta) = |K\vert\theta^{\perp}|.
\end{equation}
If $L$ is the projection body of some convex body, we simply say 
that $L$ is a projection body.  The Minkowski (vector) sum of projection bodies
is also a projection body. Every projection body is the limit in the Hausdorff metric
of Minkowski sums of symmetric intervals. An origin-symmetric convex body
in $\R^n$ is a projection body if and only if the polar body is the unit ball of
an $n$-dimensional subspace of $L_1;$ see \cite{S2,G3,K4} for proofs and 
more properties of projection bodies.
\medbreak
\noindent {\bf Proof of Theorem \ref{main-proj1}.} By approximation (see \cite[Theorem 3.3.1]{S2}), 
we can assume that $K,L$ 
are infinitely smooth. It was proved in \cite{KRZ} that
\begin{equation} \label{f-proj}
P_K(\xi)=|K\vert \xi^\bot| = -\frac 1{\pi} \widehat{f_K}(\xi),\qquad \forall \xi\in S^{n-1},
\end{equation}
where $f_K$ is extended from the sphere to a homogeneous function of degree 
$-n-1$ on the whole $\R^n,$ and the Fourier transform $\widehat{f_K}$ is the 
extension of a continuous function $P_K$ on the sphere to a homogeneous of degree 1
function on $\R^n.$

Therefore, the condition (\ref{proj2}) can be written as
\begin{equation} \label{fourier-proj}
\frac 1{\pi} \widehat{f_K}(\xi) \ge  \frac 1{\pi} \widehat{f_L}(\xi) + \e, \qquad \forall \xi\in S^{n-1}.
\end{equation}
It was also proved in \cite{KRZ} that an infinitely smooth origin-symmetric convex body 
$L$ in $\R^n$ is a projection body if and only if 
$\widehat{h_L} \le 0$ on the sphere $S^{n-1}.$ Therefore, integrating (\ref{fourier-proj})
with respect to a negative density,
$$\int_{S^{n-1}} \widehat{h_L}(\xi) \widehat{f_L}(\xi)\ d\xi \ge \int_{S^{n-1}} \widehat{h_L}(\xi) \widehat{f_K}(\xi)\ d\xi 
- \pi\e \int_{S^{n-1}} \widehat{h_L}(\xi)\ d\xi.$$
Using this, (\ref{vol-proj}) and (\ref{pars-proj}), we get
$$ (2\pi)^n n |L| = (2\pi)^n \int_{S^{n-1}} h_L(x) f_L(x)\  dx =
\int_{S^{n-1}} \widehat{h_L}(\xi) \widehat{f_L}(\xi)\ d\xi$$
$$\ge \int_{S^{n-1}} \widehat{h_L}(\xi) \widehat{f_K}(\xi)\ d\xi  - \pi\e\int_{S^{n-1}}\widehat{h_L}(\xi)\ d\xi$$
\begin{equation} \label{eq31}
=(2\pi)^n \int_{S^{n-1}} h_L(x) f_K(x)\  dx - \pi\e\int_{S^{n-1}}\widehat{h_L}(\xi)\ d\xi.
\end{equation}
We estimate the first summand from below using the first Minkowski inequality:
\begin{equation} \label{eq32}
(2\pi)^n \int_{S^{n-1}} h_L(x) f_K(x)\  dx \ge (2\pi)^n n \left(\vol_n(L)\right)^{\frac 1n} \left(\vol_n(K)\right)^{\frac {n-1}n}.
\end{equation}

To estimate the second term in (\ref{eq31}), note that, by (\ref{f-proj}),
the Fourier transform of the curvature function of the Euclidean ball
$$\widehat{f_2}(\xi) = -\pi |B_2^{n-1}|.$$
Therefore, by Parseval's formula, (\ref{mixvol-proj})  and the first Minkowski inequality,
$$\pi \e \int_{S^{n-1}}\widehat{h_L}(\xi)\ d\xi = - \frac {\e}{|B_2^{n-1}|}
\int_{S^{n-1}}\widehat{h_L}(\xi)\widehat{f_2}(\xi)\ d\xi$$
$$= -  \frac {(2\pi)^n \e}{|B_2^{n-1}|}
\int_{S^{n-1}} h_L(x) f_2(x)\ dx
= -  \frac {(2\pi)^n \e}{|B_2^{n-1}|} nV_1(B_2^n,L)$$
$$\le  - \frac {(2\pi)^n n\e}{|B_2^{n-1}|} |L|^{\frac 1n}
 |B_2^n|^{\frac {n-1}n}
 = - (2\pi)^n n \e c_n |L|^{\frac 1n}.$$
 Combining this with (\ref{eq31}) and (\ref{eq32}), we get the result.
 \endpf
 
As explained in the Introduction, the separation result of Theorem \ref{main-proj1} leads to
a volume difference inequality of the type (\ref{diff2}). 

\begin{co} If $L$ is a projection body in $\R^n$ and $K$ is an arbitrary origin-symmeric convex body in $\R^n$
so that $$\min_{\xi\in S^{n-1}} (|L\vert \xi^\bot| - |K\vert \xi^\bot|) >0,$$ then
\begin{equation} \label{diff-proj}
|L|^{\frac {n-1}n} - |K|^{\frac {n-1}n} \ge c_n \min_{\xi\in S^{n-1}} (|L\vert \xi^\bot| - |K\vert \xi^\bot|).
\end{equation}
\end{co}
Putting $K=\beta B_2^n$ in (\ref{diff-proj}) and sending $\beta \to 0,$ we get a hyperplane inequality of the type (\ref{hyper2}),
which was earlier deduced directly from the solution to Shephard's problem in \cite[Corollary 9.3.4]{G3}:
if $L$ is a projection body in $\R^n$, then
\begin{equation} \label{hyper-proj}
|L|^{\frac {n-1}n} \ge c_n \min_{\xi\in S^{n-1}} |L\vert \xi^\bot|.
\end{equation}
Recall that $c_n > 1/\sqrt{e}.$ For general symmetric convex bodies, Ball \cite{Ba2} proved that 
$c_n$ may and has to be replaced in (\ref{hyper-proj}) by
 $c/\sqrt{n},$ where $c$ is an absolute constant. 
Also, note that the inequality
$$|L|^{\frac {n-1}n} \le c_n \max_{\xi\in S^{n-1}} |L\vert \xi^\bot|$$
holds for all origin-symmetric convex bodies and follows from the Cauchy projection formula for the surface area 
(the first part of (\ref{cauchy})) and the classical isoperimetric inequality; see \cite[p. 363]{G3}.

The volume difference inequality (\ref{diff-proj}) allows to prove a hyperplane
inequality for the surface area of projection bodies.
\begin{co} \label{surf-proj}Let $L$ be a projection body in $\R^n,\ n\ge 3,$ then the surface area
$$S(L) \ge \frac n{n-1} c_n \min_{\xi\in S^{n-1}} S(L\vert\xi^\bot)\ |L|^{\frac 1n}.$$
\end{co}
\pf The surface area of $L$ can be computed as
$$S(L) = \lim_{\e \to +0} \frac {\left|L+\e B_2^n\right| - \left|L\right|}{\e}.$$
For every $\e>0$ the Minkowski sum $L+\e B_2^n$ is also a projection body. 
The inequality (\ref{diff-proj}) with the bodies
$L+\e B_2^n$ and $L$ in place of $L$ and $K$ implies
\begin{equation} \label{surflimit}
\frac {|L+\e B_2^n|^{\frac {n-1}n} - |L|^{\frac {n-1}n}}{\e} \ge
c_n \min_{\xi\in S^{n-1}} \frac {|(L\vert \xi^\bot)+\e B_2^{n-1}| - 
|L\vert \xi^\bot|}{\e}.
\end{equation}
By the Minkowski theorem on mixed volumes (\cite[Theorem 5.1.6]{S2} or \cite[Theorem A.3.1]{G3}),
\begin{equation} \label{quer}
\frac{|(L\vert \xi^\bot)+ \e B_2^{n-1}|-|L\vert \xi^\bot|}{\e} = 
\sum_{i=1}^{n-1} {n-1 \choose i} W_i(L\vert \xi^\bot) \e^{i-1},
\end{equation}
where $W_i$ are quermassintegrals. The function $\xi\mapsto L\vert \xi^\bot$ is continuous 
from $S^{n-1}$ to the class of origin-symmetric convex sets equiped with the Hausdorff metric,
and $W_i$'s are also continuous with respect to this metric (see \cite[p. 275]{S2}), so the 
functions $\xi\mapsto W_i(L\vert \xi^\bot)$ are continuous and, hence, bounded on the sphere.
This implies that the left-hand side of (\ref{quer}) converges to $S(L\vert \xi^\bot),$ as $\e\to 0,$
uniformly with respect to $\xi.$ The latter allows to switch the limit and maximum  
in the right-hand side of (\ref{surflimit}), as $\e\to 0$.
Sending $\e$ to zero in (\ref{surflimit}), we get
$$\frac {n-1}n |L|^{-1/n} S(L) \ge c_n  \min_{\xi \in S^{n-1}} S(L\vert \xi^\bot).$$
\endpf
Note the similarity between Corollaries \ref{surf-proj} and \ref{hyper-aversect}. In fact, the Cauchy projection formula
(see for example \cite[p. 408]{G3}) can be written as
\begin{equation}\label{cauchy}
S(L) = \frac 1{|B_2^{n-1}|} \int_{S^{n-1}} |L\vert \xi^\bot| d\xi = \frac {|S^{n-1}|}{|B_2^{n-1}|}\  {\rm{ap}}(L),
\end{equation}
where we denote by ${\rm{ap}}(L)$ the average hyperplane projection of $L.$
Thus, the inequality of Corollary \ref{surf-proj} turns into
$${\rm{ap}}(L) \ge \frac {|B_2^{n-1}|}{ |B_2^{n-2}|  |B_2^n|^{\frac 1n}}  \min_{\xi\in S^{n-1}} {\rm{ap}}(L\vert \xi^\bot) \ |L|^{\frac 1n}.$$
\bigbreak

A stability result for hyperplane projections was proved in \cite[Theorem 5]{K6}. Define the normalized circumradius of $L$ by
$$R(L) = \frac{\max_{\xi\in S^{n-1}} \rho_L(\xi)}{|L|^{\frac 1n}}.$$

\begin{theorem} \label{main-proj} (\cite{K6}) Suppose that $\e>0$,  $K$ and $L$ are origin-symmetric
convex bodies in $\R^n,$ and $L$ is a projection body.  If for every $\xi\in S^{n-1}$
\begin{eqnarray}\label{proj1}
|K\vert \xi^\bot| \le |L\vert \xi^\bot| +\e,
\end{eqnarray}
then
$$|K|^{\frac{n-1}n}  \le |L|^{\frac{n-1}n} + \sqrt{\frac{2\pi}n}\ R(L) \e.$$
\end{theorem}

Since the answer to Shephard's problem is negative in most dimensions, one can ask what
condition on the hyperplane projection function does imply the inequality for volumes.
Yaskin \cite{Y} proved that for $\alpha \in [n,n+1)$ the inequalities
\begin{eqnarray}\label{eqn:condition2}
(-\Delta)^{\alpha/2} P_K(\xi)\ge (-\Delta)^{\alpha/2} P_L(\xi), \qquad \forall \xi\in S^{n-1}
\end{eqnarray}
imply that $|K| \le |L|,$ where the projection functions are extended to homogeneous functions 
of degree 1 on the whole $\R^n.$ The latter result is no longer true for $\alpha < n.$
We end this section by formulating the stability version of the result of Yaskin.
 
 \begin{theorem} \label{main-Y} Let $\e>0,\ \alpha \in [n,n+1)$,  $K$ and $L$ be origin-symmetric
infinitely smooth convex bodies in $\R^n$, $n\ge 3$, so that for every $\xi\in S^{n-1}$
$$(-\Delta)^{\alpha/2} P_K(\xi)\le (-\Delta)^{\alpha/2} P_L(\xi) +\e.$$
Then
$$|K|^{\frac{n-1}n}  \le |L|^{\frac{n-1}n} + c \e,$$
where 
$$c=\frac{\Gamma(\frac{n-\alpha+1}2) \left|S^{n-1}\right| R(L)}{2^{\alpha+1}\pi^{\frac n2}
\Gamma(\frac{\alpha+1}2)n}.$$
\end{theorem}
Note that this is no longer true if $\alpha<n,$ because the underlying comparison result fails, 
as shown in \cite{Y}.

\section{Arbitrary measures} Zvavitch \cite{Zv} found a remarkable generalization of the Busemann-Petty problem
to arbitrary measures, namely, one can replace volume by any measure with 
even continuous density in $\R^n.$ In particular, if $K$ is an intersection body in $\R^n$  and 
$L$ is an arbitrary origin-symmetric star body in $\R^n,$ then the inequalities
$$\mu(K\cap \xi^\bot) \le \mu(L\cap \xi^\bot), \qquad \forall \xi\in S^{n-1}$$
imply
$$\mu(K)\le \mu(L). $$

Stability in Zvavitch result was established in \cite[Theorem 2]{K8}. Note that in the case of volume (when $f \equiv 1$),
the result of Theorem \ref{stab-measure} is weaker than that of Theorem \ref{main-int}. Also, Theorem \ref{stab-measure}
was formulated in \cite{K8} for dimensions up to 4 only, however, the proof works in any dimension under the 
assumption that $K$ is an intersection body.
\begin{theorem} \label{stab-measure} (\cite{K8}) Let  $f$ be an even non-negative continuous function on $\R^n,$
let $\mu$ be the measure with density $f,$ let $K$ and $L$ be origin-symmetric star
bodies in $\R^n,$ and let $\e>0.$ Suppose that $K$ is an intersection body and that for every $\xi\in S^{n-1},$
\begin{equation} \label{cond}
\mu(K\cap \xi^\bot) \le \mu(L\cap \xi^\bot) +\e.
\end{equation}
Then
\begin{equation}\label{concl}
\mu(K)\le \mu(L) + \frac {n}{n-1}c_n |K|^{1/n} \e.
\end{equation}
\end{theorem}

Interchanging $K$ and $L,$ we get the volume difference inequality.
\begin{co} \label{ineq-meas} If $K$ and $L$ are intersection bodies in $\R^n$ (in particular, 
any origin-symmetric convex bodies in $\R^n,\ n\le 4$), then 
$$\left|\mu(K) - \mu(L)\right| $$
\begin{equation} \label{measure}
\le \frac {nc_n}{n-1}\max_{\xi \in S^{n-1}}\left|\mu(K\cap \xi^\bot) - \mu(L\cap \xi^\bot)\right| 
 \max\left\{|K|^{\frac 1n}, |L|^{\frac1n}\right\}.
\end{equation}
\end{co}
Sending $L$ to the emply set, we arrive at the hyperplane inequality for arbitrary measure.
\begin{co}\label{main-measure} If $K$  is an intersection body in $\R^n$ (in particular, any origin-symmetric convex 
body in $\R^n,\ n\le 4$), then 
\begin{equation} \label{arbmeas}
\mu(K) \le \frac n{n-1} c_n \max_{\xi \in S^{n-1}} \mu(K\cap \xi^\bot)\ |K|^{1/n}.
\end{equation}
\end{co}

The constant in (\ref{arbmeas}) is sharp, it is achieved asymptotically when $K=B_2^n$
and $\mu$ converges weakly to the uniform measure on the sphere $S^{n-1};$  see \cite{K8}.

\bigbreak

{\bf Acknowledgements.} This work was partially supported by
the US National Science Foundation through 
grant DMS-1001234. I wish to thank the Max Planck
Institute for Mathematics for support and hospitality during my
stay in Spring 2011, when this project was initiated.

\end{document}